\documentclass[twoside,11pt]{article}

\usepackage{float}

\usepackage[utf8]{inputenc}
\usepackage{graphicx, subfigure}
\usepackage[top=1in,bottom=1in,left=1in,right=1in,footnotesep=0.5in]{geometry}
\usepackage[sort&compress]{natbib} 
\usepackage{amsfonts, amsmath, amssymb, amsthm, constants, bbm}
\usepackage{MnSymbol}
\usepackage{booktabs} 
\usepackage{enumitem}

\usepackage{dsfont} 

\usepackage[utf8]{inputenc}
\usepackage{MnSymbol}
\usepackage{booktabs} 
\usepackage{enumitem}
\usepackage{amsfonts, amsmath, amsthm}

\usepackage{color}
\definecolor{darkred}{RGB}{150,0,0}
\definecolor{darkgreen}{RGB}{0,150,0}
\definecolor{darkblue}{RGB}{0,0,150}

\usepackage{hyperref}
\hypersetup{colorlinks=true, linkcolor=darkred, citecolor=darkgreen, urlcolor=darkblue}
\usepackage{url}

\usepackage[sort&compress]{natbib}

\theoremstyle{remark}

\def\beq{\begin{equation}} 
\def\eeq{\end{equation}}
\def\beqn{\begin{eqnarray*}}
\def\eeqn{\end{eqnarray*}}
\def\Bitem{\begin{itemize}\setlength{\itemsep}{.2in}}
\def\bitem{\begin{itemize}\setlength{\itemsep}{.05in}}
\def\eitem{\end{itemize}}
\def\Benum{\begin{enumerate}\setlength{\itemsep}{.2in}}
\def\benum{\begin{enumerate}\setlength{\itemsep}{.05in}}
\def\eenum{\end{enumerate}}
\def\bmult{\begin{multline*}}
\def\emult{\end{multline*}}
\def\bcenter{\begin{center}}
\def\ecenter{\end{center}}
\def\bframe{\begin{frame}}
\def\eframe{\end{frame}}

\newcommand{\secref}[1]{Section~\ref{sec:#1}}
\newcommand{\figref}[1]{Figure~\ref{fig:#1}}
\newcommand{\tabref}[1]{Table~\ref{tab:#1}}

\DeclareMathOperator*{\argmax}{arg\, max}

\newcommand{\blambda}{{\boldsymbol\lambda}}

\newcommand{\btheta}{{\boldsymbol\theta}}

\def\bbR{\mathbb{R}}

\def\1{\mathbbm{1}}
\newcommand{\IND}[1]{\1_{\{ #1 \}}}

{
  \theoremstyle{plain}
  
}
{
  \theoremstyle{plain}
  
}
{
  \theoremstyle{plain}
  
}
{
  \theoremstyle{remark}
  \newtheorem{remark}{Remark}
}

\begin{document}

\thispagestyle{empty}
\title{Fitting a Multi-modal Density by Dynamic Programming}
\author{
Ery Arias-Castro\,\footnote{Department of Mathematics, University of California, San Diego, USA \newline \indent \quad \url{https://math.ucsd.edu/\~eariasca}} 
\and 
He Jiang\,\footnote{Department of Mathematics and Statistics, California State Polytechnic University, Pomona, USA}
}
\date{}
\maketitle

\begin{abstract}
We consider the problem of fitting a probability density function when it is constrained to have a given number of modal intervals. We propose a dynamic programming approach to solving this problem numerically. When this number is not known, we provide several data-driven ways for selecting it. We perform some numerical experiments to illustrate our methodology. 
\end{abstract}

\section{Introduction}  \label{sec:intro}

Density estimation is an important task in exploratory data analysis, for example, the histogram is introduced in the most basic courses in statistics, and as such has been the object of some longstanding and intense study in statistics and related fields \citep{sheather2004density, silverman2018density, scott2015multivariate}. Rather than assuming that the underlying density comes from a parametric family (e.g., normal), or some smoothness class (e.g., differentiable with some given bound on the derivative), we focus here on the shape constraint that the density has at most a given number of modal intervals.
A mode is an important feature and defining the density based on the number of modal intervals is thus intuitive. It also goes hand-in-hand with a modal approach to clustering, in particular, as proposed by \citet{fukunaga1975estimation}. 

In this paper, we propose a dynamic programming way of fitting a multi-modal density to a numerical (one-dimensional) sample, and discuss data-driven ways of choosing the number of modal intervals when this information is unknown.

\subsection{Importance of density modes}
\label{sec:importance_modes_subsec}

\begin{quote}
{\em Thus the `mean,' the `mode,' and the `median' have all distinct characters important to the statistician.}
\hfill \cite{pearson1895contributions}
\end{quote}

\begin{figure}[!htpb]
\centering
\includegraphics[scale=0.3]{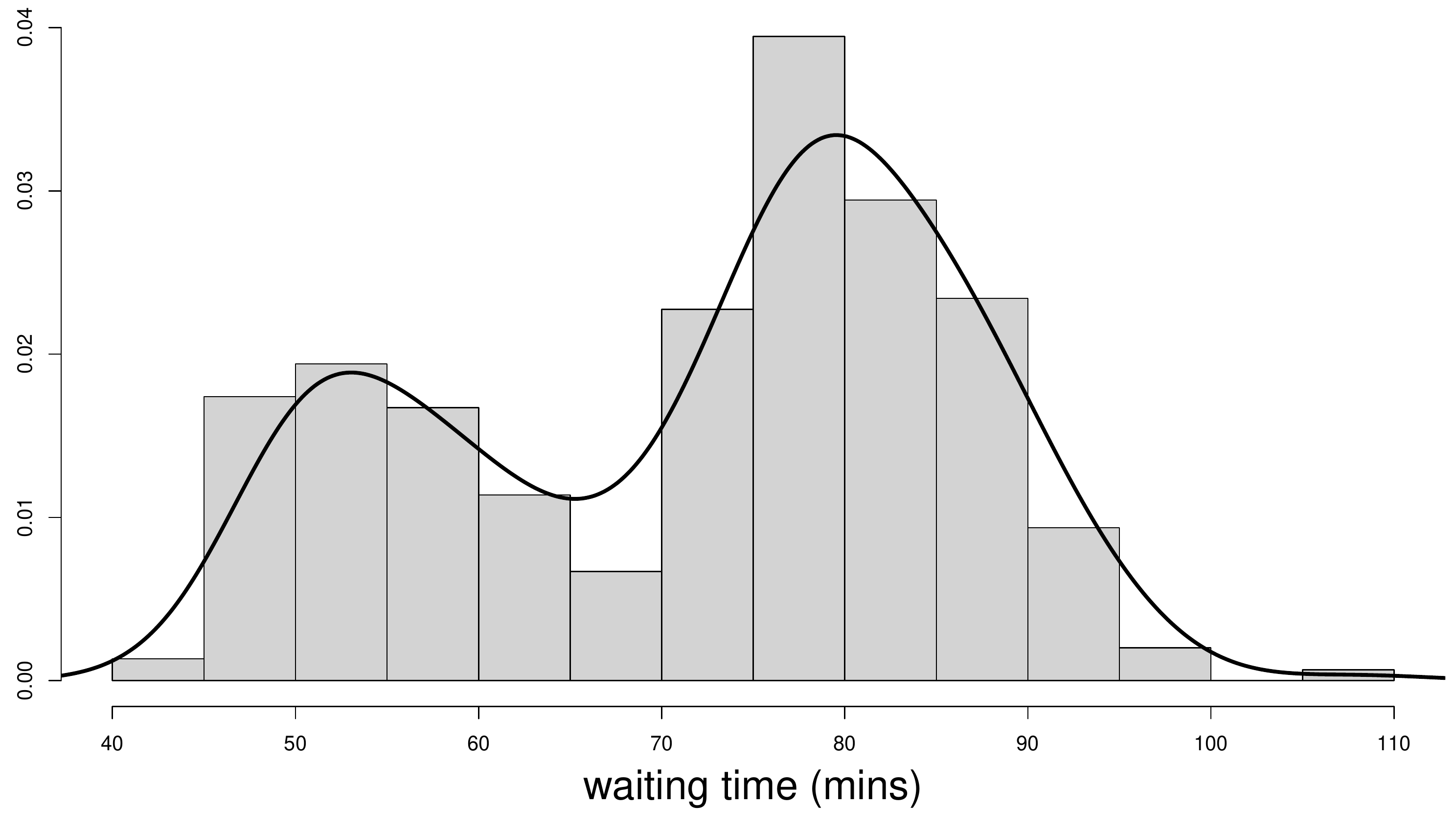}
\caption[An example in the Geyser dataset.]{A histogram of the waiting times between eruptions of the Old Faithful Geyser, Yellowstone National Park in Wyoming, United States, with an overlay of a kernel density estimate. (Based on the {\sf geyser} dataset in the {\sf MASS} package in {\sf R}.)}
\label{fig:geyser_histogram_density}
\end{figure}

Modes are important features of densities, and do not fail to catch the eye of the analyst when they manifest themselves in plots of histograms or other estimates of the density. See \figref{geyser_histogram_density} for a classical example.
Estimating the modes of the underlying density has thus been an important problem in statistics dating back to \cite{parzen1962estimation}, who proposed to first estimate the density by kernel density estimation --- pioneered just a few years by by \cite{rosenblatt1956remarks} --- and then locate the modes of the resulting estimator. 
The problem of estimating the density modes has been the subject of great, ongoing interest by statisticians since then; see \citep{chacon2020modal} for a recent review.

In addition to being remarkable features of densities, modes are also useful when the goal is clustering. Working at the population level as in \citep{chacon2015population}, a clustering amounts to a partitioning of the population, and \citet{fukunaga1975estimation} proposed to partition the population according to the basins of attraction of the gradient ascent flow defined by the density. These are called the ``descending manifolds'' in Morse theory \citep{chen2017statistical}, and the ones that pertain to modes provide a partition of the support of the density up to a set of measure zero when the density is regular in a certain sense (specifically, when it is of Morse type). The corresponding methodology is known as the ``mean-shift algorithm'' \citep{cheng1995} and is now a well-established approach to clustering with its own literature \citep{carreira2015review, menardi2016review, chacon2020modal}.
This approach to clustering turns out to be very closely related to the level set/cluster tree approach of \cite{hartigan1975clustering}, which itself has generated a good amount of research \citep{stuetzle2010generalized, rinaldo2012stability, chaudhuri2014consistent}.
For a connection between these two approaches to clustering, see our recent work \citep{arias2021asymptotic, arias2021moving}. The level set approach to clustering is also known to be intimately connected to single-linkage clustering \citep{hartigan1981consistency, penrose1995single}.
We note that, in dimension one --- our focus here --- modal clustering amounts to partitioning the support of the density according to the intervals defined by local minima (indeed, not the modes or local maxima), and thus is exceedingly simple to execute. In practice, of course, inference is based on a sample, and these local minima need to be estimated. 
When fitting a multi-modal density, estimates for these local minima are automatically available.

Density modes are also of interest in mixture models. For example, the number of modes in a Gaussian mixture model has been the subject of some attention \citep{mclachlan2014number}. We note, however, that the consideration of modes is not particularly natural in the context of mixtures, for the simple reason that the modes are not directly related to the component centers.
Also, when components are too close to each other in a mixture model, it is difficult to interpret the meaning of each component. 
When fitting a multi-modal density, the interpretation is arguably clearer.

\subsection{Estimating a uni-modal or multi-modal density}

As we just alluded to, a model for multi-modal densities offers a competitive alternative to a mixture model. In this paper, we simply choose as model the entire class of densities with at most $K$ modes. While a {\em mode} typically means an isolated local maximum, to better accommodate step density functions --- which are important in many respects, including in their connection to maximum likelihood --- we say that a density $f$ is {\em $K$-modal} if there are 
$\lambda_1 < \cdots < \lambda_{K+1}$ such that $f$ is uni-modal on $[\lambda_{k}, \lambda_{k+1}]$ for all $k=1,\dots,K$, where $\lambda_1 = -\infty$ and $\lambda_{K+1} = \infty$.
We refer to $\lambda_{2}, \lambda_{3}, \dots, \lambda_{K}$ as {\em knots}. We note that they are local minima, but also more than that.
To be sure, $f$ is {\em uni-modal} on some interval $(a,b)$ if there is some $a \le \mu \le b$ such that $f$ is non-decreasing on $(a,\mu]$ and non-increasing on $[\mu, b)$. In that case, $\mu$ is a mode of the density $f$.
We call an interval where $f$ is uni-modal a {\em modal interval}.

\begin{remark}
As defined, the modal intervals may not be unique. In fact, the number of modal intervals, $K$, is also not uniquely defined. However, we can remedy this by letting $K$ be the smallest integer such that there is a partition of the real line into $K$ modal intervals for $f$. (Note that $K=\infty$ is possible.) And then, when $K < \infty$ and $\lambda_{2}, \lambda_{3}, \dots, \lambda_{K}$ are knots, we may redefine $\lambda_k$ by taking the midpoint of the interval containing $\lambda_k$ where $f$ takes the value $f(\lambda_k)$.
\end{remark}

Given the significance of modes, there is some amount of literature on fitting densities with some constraints on the number of modes, although most of the attention has gone to uni-modal densities, i.e., densities with a single mode, with both methodological and theoretical developments. This effort dates back to the work of \cite{grenander1956theory} on the estimation of a monotonic density, and his derivation of the corresponding maximum likelihood estimator, which nowadays bears his name. A monotone density is, of course, uni-modal. Vice versa, a uni-modal density with known mode at $\mu$ is non-decreasing on $(-\infty, \mu]$ and non-increasing on $[\mu, \infty)$, so that the maximum likelihood is given by Grenander's estimator applied on each side of $\mu$ with the proper monotonicity constraint. 
\cite{rao1969estimation} obtained the asymptotic distribution of this estimator pointwise when applied to estimating uni-modal densities with a known mode, and showed this estimator's consistency on intervals not containing the mode.
When the mode is unknown, the likelihood is simply maximized over the location $\mu \in \bbR$.
As an extension, \cite{wegman1970maximum1, wegman1970maximum2} proposed a maximum likelihood estimator of a uni-modal density when the location of the mode is unknown, and showed that for sufficiently large sample size and on certain regions, this estimator agrees with the maximum likelihood estimator of a uni-modal density with known mode. 
We note that penalized versions of Grenander's estimator have also been proposed \citep{woodroofe1993penalized,sun1996adaptive}, and their use could be extended to fitting a uni-modal density.
Besides maximum likelihood methods, the problem of uni-modal density estimation has also been considered using other approaches, including via the use of splines \citep{bickel1996some, meyer2004consistent, meyer2012nonparametric}, data sharpening \citep{choi1999miscellanea, braun2001data,wolters2012greedy}, and optimal transport \citep{cumings2018shape}.
We also mention \citep{birge1987estimating, birge1987risk}, where nonparametric minimax rates are derived for the problem of estimating a uni-modal density or decreasing density using well chosen histograms with unequal bin width; and \citep{hengartner1995finite}, where conservative finite-sample confidence regions are computed for the entire density under monotonicity or uni-modality. 

Given the amount of literature on fitting a uni-modal density, it is natural to consider fitting multi-modal density with any given number of modal intervals. There is a significant literature on testing for multi-modality, which we review later in the paper in \secref{testing multi-modality}. In terms of estimation methodology, the literature is much more scarce.
\cite{minnotte1993mode} proposed an exploratory, multi-scale method they named the {\em mode tree}. A similar method, called {\em SiZer}, was later proposed by \cite{chaudhuri1999sizer, chaudhuri2000scale}, and later, a rank-based variant was suggested by \cite{dumbgen2008multiscale}. Such methods, although they do not offer a way of estimating a multi-modal density, may be seen as offering a way to partition the real-line into (possible) modal intervals.
Directly tackling the estimation problem, \cite{wolters2018enforcing} considered fitting densities with $1$ or $2$ modes by adding a function to an initial kernel density estimate and encouraging the resulting function to satisfy the desired constraint using a quadratic program. \cite{dasgupta2021modality} considered fitting multi-modal densities by first fitting a template density with the given number of modes, in the spirit of \citep{cheng1999nonparametric}, and then applying a shape-preserving transformation to obtain the final estimate. 

\subsection{Contribution and content}

We consider fitting a probability density function with a constraint on the number of modal intervals. We introduce a dynamic programming approach to solving this problem numerically which, in principle, may be based on any method for fitting a uni-modal density on a given interval. (In our numerical experiments we use the method recently proposed by \cite{wolters2018enforcing}.) To lighten up the computational burden, we propose a multi-grid search approach. 
In addition, we offer data-driven ways for selecting the number of modal intervals when this information is unknown.


The organization of the paper will be as follows. In \secref{dynamic_programming}, we describe our method for fitting a $K$-modal density via dynamic programming, where $K$ is provided by the user. We also introduce a multi-grid search approach for improving the estimate of the knots. In \secref{choosing_number}, we provide several data-driven options for choosing the number of modal intervals $K$, when this information is unknown. In \secref{k_mode_experiment}, we report on some numerical experiments. We end the paper with a discussion in \secref{conclusion_kmode}.

\section{Dynamic programming method}
\label{sec:dynamic_programming}

We have at our disposal a sample of real-valued observations.
These data points are assumed ordered (without loss of generality).  They are denoted $x_1 \le \cdots \le x_n$ and gathered in a vector $\mathbf{x} = (x_1, \dots, x_n) \in \mathbb{R}^n$.
Given $\mathbf{x}$, our goal is to fit a $K$-modal density. 

\subsection{Maximum likelihood estimation}

It is natural to first consider an estimation by maximum likelihood, which we can formalize as the following optimization problem:
\begin{equation}
\label{original_optimization}
\text{maximize} \ \sum_{i=1}^n \log (f(x_i)),
\end{equation}
over $f$ an arbitrary $K$-modal density.

While when $K=1$ this leads to the Grenander estimator for fitting a uni-modal density, when $K \ge 2$, this problem is not well-posed due to the fact that the likelihood is not bounded over the class of $K$-modal densities. In fact, this is already the case when fitting mixtures of two or more components, including Gaussian mixtures, if the variances are left unconstrained. See, e.g., \citep[Sec 7]{day1969estimating} or \citep[Ex 5.51]{van2000asymptotic}. 
We note that, for fitting Gaussian mixtures, some solutions to this issue have been suggested. For example, \cite{hathaway1986constrained} proposed a constrained EM algorithm when maximizing the likelihood with a lower bound on component weights and an upper bound on the aspect ratio of the component variances, which allows the maximum likelihood estimator to be consistent \citep{hathaway1985constrained}.

In our situation, we can think of multiple ways to regularize the optimization problem. One way to do so is to constrain the search for modes to a finite grid of values that avoids all the data points.
A further discretization of the problem would involve a search over a finite grid (the same one, or another one) for the $K-1$ knots defining the $K$ candidate modal intervals. For such a $(K-1)$-tuple, the optimization would lead to fitting a uni-modal density on each of the corresponding intervals with a search for the mode done over the mode grid. 

\subsection{Fitting uni-modal densities}
\label{sec:fitting_unimodal}

\begin{figure}[htpb]
\centering
\includegraphics[scale=0.5]{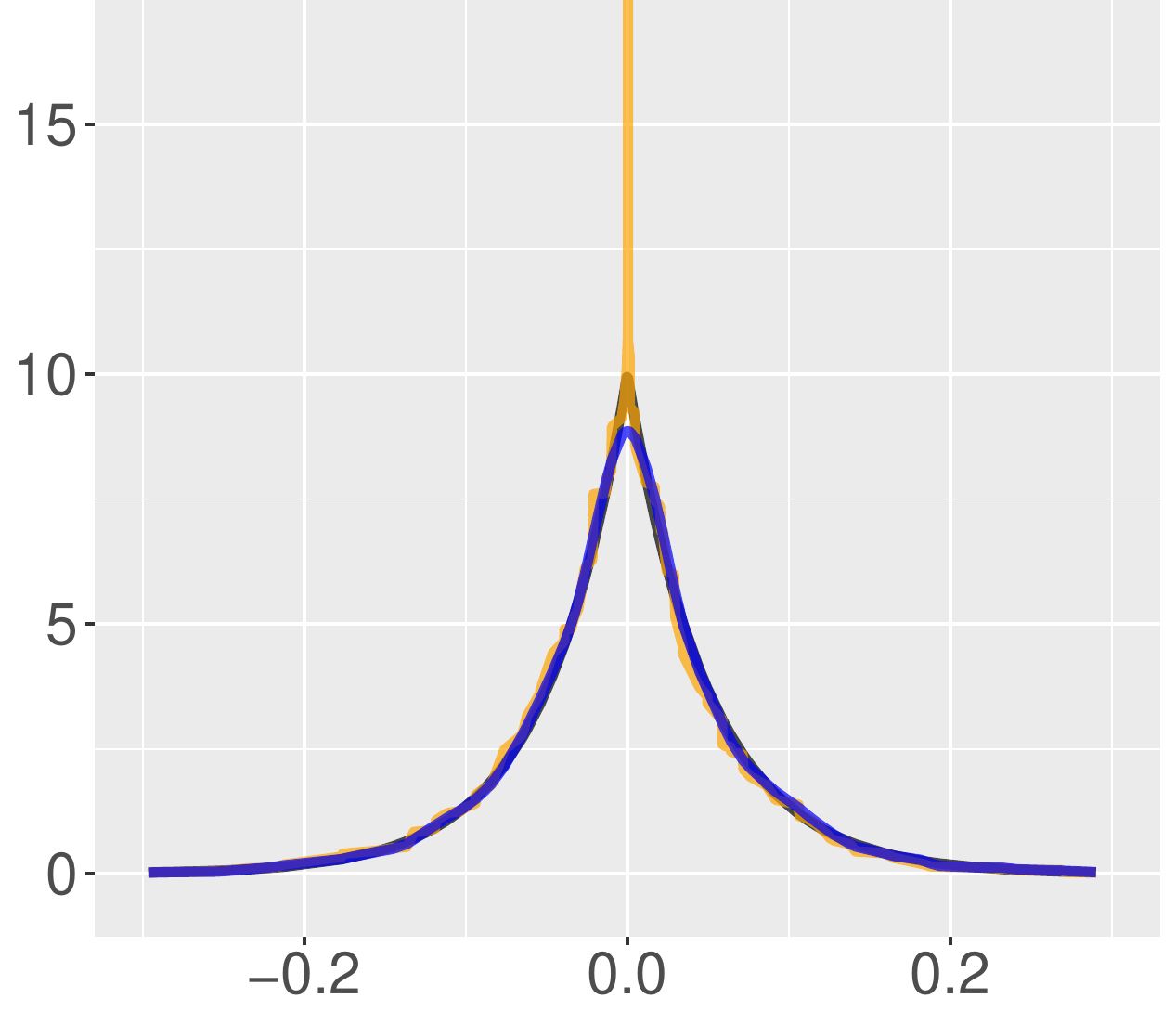}
\caption[Uni-modal densities by maximum likelihood and adjusted KDE.]{ The density $f$ of the Laplace distribution with  mean $0$ and standard deviation $0.05$, in black, overlaid with the uni-modal density acquired via maximum likelihood \citep{grenander1956theory}, in orange, and the uni-modal density acquired via adjusted KDE \citep{wolters2018enforcing}, in blue. The plots does not show the entire MLE estimator, as it takes a value larger than $3\times 10^4$ at the mode.}
\label{fig:k_modes_ml_high_density}
\end{figure}

Even though the maximum likelihood estimator --- at least in its constrained and discretized form --- exists, it is known to have an issue, which comes from the Grenander estimator being inconsistent at the mode \citep{balabdaoui2011grenander, woodroofe1993penalized, bickel1996some}. Although it is consistent everywhere else and so,  practically speaking, this is not particularly important, the resulting fit  is not ``visually'' pleasing, as the value at the mode tends to be much larger than the other values that the density estimator takes, resulting in a ``spike'' at the mode. See \figref{k_modes_ml_high_density} for a numerical illustration.
As our approach can accommodate any method for fitting a uni-modal density, in our implementation we chose another method, that of \cite{wolters2018enforcing}.  


\subsection{Problem decomposition}
\label{sec:decomposition_uni-modal}

We approach the optimization problem \eqref{original_optimization} in the following way, which involves a regularization by discretization and the use of a well-behaved method for fitting a uni-modal density on a given sample. For a candidate density $f$, let $\blambda = ( \lambda_1, \dots, \lambda_{K+1} )$ with $\lambda_1 = -\infty$ and $\lambda_{K+1} = \infty$, and with $\lambda_2 < \cdots < \lambda_{K}$ denoting some knots defining modal intervals for $f$.

\begin{remark}
\label{rem:one_mode}
Here and afterwards, if $K=1$, there is no need to search for the knot locations, as in that case $\lambda_1 = -\infty$ and $\lambda_{2} = \infty$. We therefore proceed assuming that $K \ge 2$.
\end{remark}

Define $\pi_k = \int_{\lambda_{k}}^{\lambda_{k+1}} f(x) dx$
and assume that $\pi_k > 0$. Let $h_k$ denote the density $f$ conditional on
$[\lambda_{k}, \lambda_{k+1}]$,
i.e.,  
\begin{equation*}
h_k(x) = \frac{f(x)}{\pi_k} \IND{\lambda_{k} \leq x < \lambda_{k+1}},
\end{equation*}
and note that it is uni-modal by construction. We will refer to $\pi_k h_k$ as the `rescaled' version of $h_k$.

Realizing that the set of knots $(\lambda_2, \dots, \lambda_{K})$, the set of weights $(\pi_1, \dots, \pi_K)$, and the set of uni-modal densities $h_1, \dots, h_K$, can `move' independently of each other, we may re-express the maximization problem \eqref{original_optimization} as follows:
\begin{equation}
\label{separable_maximization_monotone}
\text{maximize} \ \sum_{k=1}^{K} \sum_{i=1}^n \log \big( \pi_{k} h_{k} (x_i) \big)  \IND{\lambda_{k} \leq x_i < \lambda_{k+1}},
\end{equation}
over $\lambda_2 < \cdots < \lambda_{K}$ chosen among grid points $g_2, \dots, g_{M}$; over $\pi_1, \dots, \pi_K > 0$ such that $\sum_{k=1}^K \pi_k = 1$; and over $h_1, \dots, h_K$ uni-modal densities with $h_k$ supported on $[\lambda_{k}, \lambda_{k+1}]$. Again, the fitting of these uni-modal densities will be eventually done by a well-behaved method, as otherwise the likelihood is not bounded and the problem above is not well-posed.

Given $\blambda$, for each $k$, let $h_k^\blambda$ denote the uni-modal density fitted by our method of choice based on the sample points falling in the interval $[\lambda_{k}, \lambda_{k+1})$. 
Note that $h_k^\blambda$ is supported on $[\lambda_{k}, \lambda_{k+1}]$.
Having computed these, the maximizing set of weights maximizes
\begin{equation}
\label{separated_pi}
\sum_{k=1}^{K} \sum_{i=1}^n \log \big( \pi_{k} h_{k}^\blambda(x_i) \big)  \IND{\lambda_{k} \leq x_i < \lambda_{k+1}}
= \sum_{k=1}^{K} n_k^\blambda \log \pi_{k} +  \sum_{k=1}^{K} \sum_{i=1}^n \log \big( h_{k}^\blambda(x_i) \big)  \IND{\lambda_{k} \leq x_i < \lambda_{k+1}},
\end{equation}
where $n_k^\blambda$ denotes the number of data points in $[\lambda_{k}, \lambda_{k+1})$. Since the second term on the right-hand side does not depend on the weights, the solution is the one maximizing the first term, which is $\pi_1^\blambda, \dots, \pi_K^\blambda$ with $\pi_k^\blambda := n_k^\blambda/n$.

The resulting estimate for the density is then 
\begin{equation*}
f^\blambda(x) = \sum_{k=1}^K \pi_k^\blambda h_k^\blambda(x) \IND{\lambda_{k} \leq x < \lambda_{k+1}}.
\end{equation*}
It thus remain to optimize with respect to the knot set $\blambda$. The search will be over a finite grid $g_2 < \cdots < g_M$, to which we add $g_1 = -\infty$ and $g_{M+1} = \infty$. 
We are able to implement this search in an efficient way.

\subsection{Dynamic programming approach}
\label{sec:dynamic_programming_description}

\begin{figure}[h]
\centering
\includegraphics[scale=0.6]{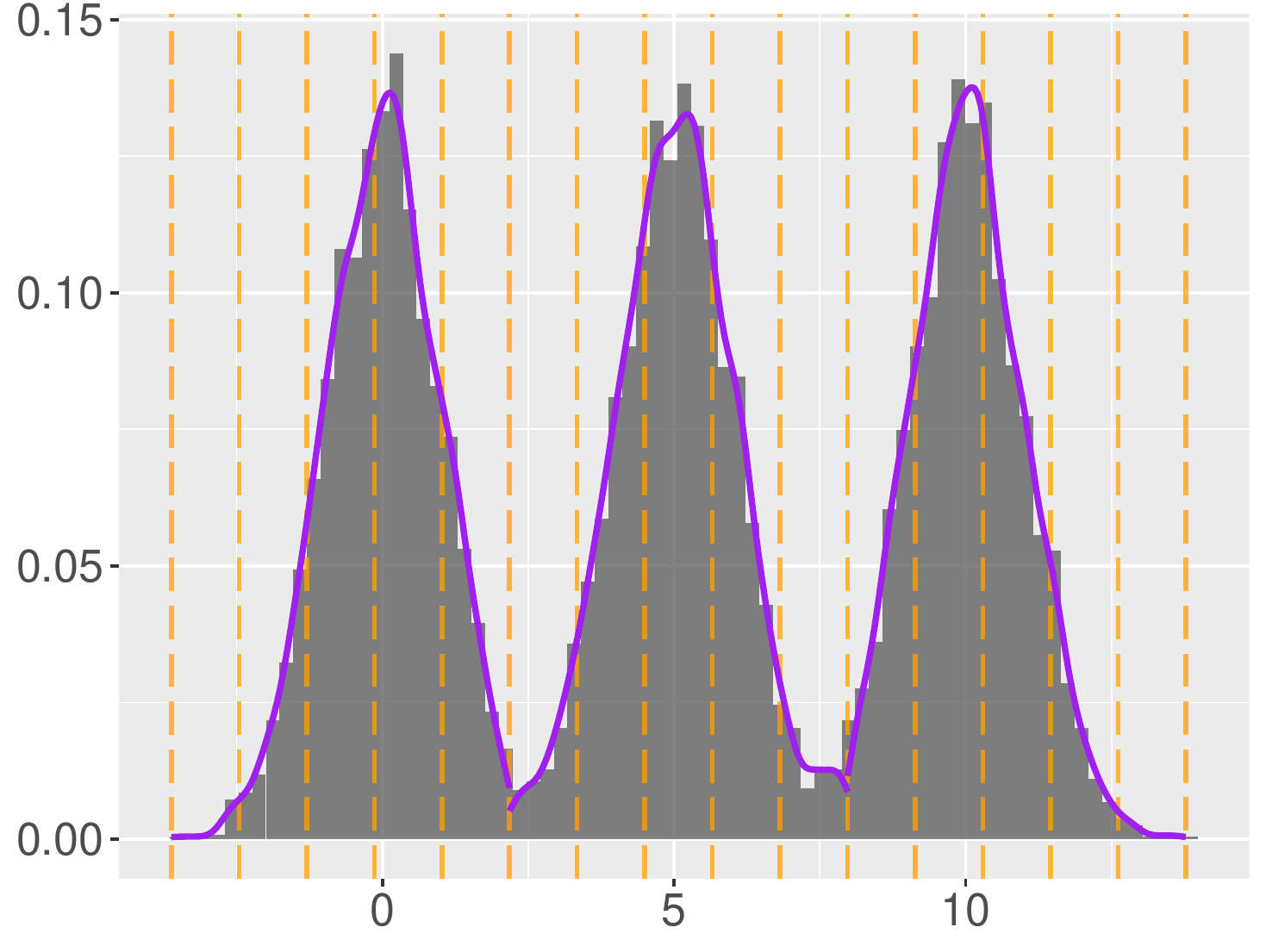}
\caption[An illustrative example on our methodology.]{ An illustrative example on our methodology. It is based on a sample of size $n = 10^4$ generated from the mixture model $\frac13\, \mathcal{N}(0,1) + \frac13 \, \mathcal{N}(5,1) + \frac13 \, \mathcal{N}(10,1)$. The number of modes for this mixture is $K=3$. In dark grey we plot a histogram of the sample. We consider the default $M = 5K - 1 = 14$ grid points as candidate locations for the knots $\blambda = (\lambda_2, \lambda_3) $. We plot these grid points, $\mathbf{g} = (g_2, \dots, g_{15})$, in vertical dotted orange lines, together with the extremes, in the left and right dotted orange lines. The $m$-th interval of interest, is located between the $m$-th and $(m+1)$-th dotted orange lines, where $m=1,2,\dots,15$. The purple line is a drawing of the density found by our dynamic programming method with knowledge of the correct number of modes. 
\label{fig:illustrative_plot}
}
\end{figure}

As described in \citep[Chap 11]{bradley1977applied}, dynamic programming is an optimization method that turns a large and complex problem into a sequence of smaller and simpler problems. In our situation, we turn a maximization problem over $K$-modal densities into a sequence of maximizations over uni-modal densities. This is possible because of the fact that in \eqref{separable_maximization_monotone} the maximization over the uni-modal densities $h_1, \dots, h_K$ can be done independently without coordination.
We provide an example to illustrate our method in action in  \figref{illustrative_plot}.

To set up the dynamic programming recursion, we define two matrices. We first let $\mathbf{S}$ be an $M$-by-$M$ upper-triangular matrix, where $\mathbf{S}[i,j]$, $i=1,\dots,M$, and $j=i,\dots,M$, denote the log-likelihood value of the (rescaled) uni-modal density fitted on $[g_{i}, g_{j+1}]$. 
In formula,
\begin{equation}
    \label{s_matrix_initialization}
    \mathbf{S}[i,j] = 
    \sum_{l=1}^n  
    \log{ \big( h_{ij}(x_l) \big) \IND{g_{i} \leq x_l < g_{j+1}}} 
    + n_{ij} 
    \log(n_{ij}/n),
\end{equation}
where $h_{ij}$ denotes the uni-modal density estimated on $[g_{i}, g_{j+1})$, and $n_{ij}$ denotes the number of data points in that interval. 

In addition to $\mathbf{S}$, let $\mathbf{D}$ be an $M$-by-$K$ lower-triangular matrix, where for $k=1,\dots, K$, $m=k,\dots, M$, $\mathbf{D}[m,k]$ denotes the maximum log-likelihood of fitting $k$ (rescaled) uni-modal densities on $[g_1, g_{m+1})$.
By definition, only entries $\mathbf{D}[m,k]$, $m \geq k$ carry meaning as it would be impossible to fit a $k$-modal density on $[g_1, g_{m+1})$ if $m < k$ as we assume that the density has at most one modal interval in any interval of the form $[g_{k}, g_{k+1})$.

We initialize $\mathbf{D}$ as follows
\begin{equation}
    \label{k_mode_d_initialization}
    \mathbf{D}[m,1] = \mathbf{S}[1,m], \ m = 1,2,\dots, M.
\end{equation}
From there, the dynamic programming recursion takes the form
\begin{equation}
\label{recursion_mode}
    \mathbf{D}[m,k] = \max_{k-1 \leq i \leq m-1} 
    \mathbf{D}[i, k-1] 
    + \mathbf{S}[i+1,m],
\end{equation}
for $k=2,\dots,K$ and $m=k,\dots, M$. 
To find the actual values of $\blambda$, we also create an $M$-by-$K$ matrix $\mathbf{I}$, where 
$\mathbf{I}[m,k]$ indicates the index of the first interval of the $k$-th uni-modal density, when fitting $k$ uni-modal densities on the first $m$ intervals. By definition, the initialization is
\begin{equation}
    \label{k_mode_i_initialization}
    \mathbf{I}[m,1] = 1.
\end{equation}
Then $\mathbf{I}[m,k]$ is updated as the 
maximizing index $i$ in \eqref{recursion_mode}. 
We provide some pseudo-code in \tabref{k_mode_dp_algorithm} and \tabref{k_mode_s_d_i_update_algorithm}.



\begin{table}[!h]
\caption{Dynamic programming approach for fitting a density with $K$ modes, which returns the maximum likelihood $\mathbf{D}[M,K]$ and the optimal knots $\blambda^\star$}
\label{tab:k_mode_dp_algorithm}
\centering
\bigskip
\setlength{\tabcolsep}{0.22in}
\begin{tabular}{ p{0.9\textwidth}  }
\toprule
{\textbf{inputs}: $\mathbf{D}$ and $\mathbf{I}$ (from \tabref{k_mode_s_d_i_update_algorithm}), number of modes $K$, grid $\mathbf{g}$ of size $M-1$} \\ 
\midrule

initialize $\mathbf{b} = \mathbf{0}_{K+1}$ and $\mathbf{b}[K+1] = M+1$

initialize $l = M$
  
\textbf{for} $k=K,\dots,1$ \textbf{do}

\hspace{3mm} 
$\mathbf{b}[k] = \mathbf{I}[l, k]$

\hspace{3mm} 
$l = \mathbf{I}[l, k] - 1$

$ {\blambda}^\star = \mathbf{g}[\mathbf{b}] $
       
\textbf{return} ${\blambda}^\star$, $\mathbf{D}[M,K]$ \\

\bottomrule
\end{tabular}
\end{table}

\begin{table}[!h]
\caption{Rountine for computing $\mathbf{D}$ and $\mathbf{I}$}
\label{tab:k_mode_s_d_i_update_algorithm}
\centering
\bigskip
\setlength{\tabcolsep}{0.22in}
\begin{tabular}{ p{0.9\textwidth}  }
\toprule
{\textbf{inputs}: data $\mathbf{x}$, number of modes $K$, grid $\mathbf{g}$ of size $M-1$
} \\ \midrule

compute $\mathbf{S}$ as in \eqref{s_matrix_initialization}

initialize $\mathbf{D}$ as in \eqref{k_mode_d_initialization} and $\mathbf{I}$ as in \eqref{k_mode_i_initialization}

\textbf{if} $K=1$ \textbf{then}

go to return step
\hspace{3mm}

\textbf{else if} $K \geq 2$ \textbf{then}

\hspace{3mm} 
\textbf{for} $k=2,\dots,K$ \textbf{do}

\hspace{3mm} \hspace{3mm} 
\textbf{for} $m=k,\dots,M$ \textbf{do}

\hspace{3mm} \hspace{3mm} \hspace{3mm} 
initialize $\mathbf{p} = \mathbf{0}_{m-k+1} $

\hspace{3mm} \hspace{3mm} \hspace{3mm} 
\textbf{for} $i=k-1,\dots,m-1$ \textbf{do}

\hspace{3mm} \hspace{3mm} \hspace{3mm} \hspace{3mm} 
$\mathbf{p}[i-k+2] = \mathbf{D}[i, k-1] + \mathbf{S}[i+1, m]$ 

\hspace{3mm} \hspace{3mm} \hspace{3mm} \hspace{3mm} 
$\mathbf{D}[m,k] = \max (\mathbf{p})$
        
\hspace{3mm} \hspace{3mm} \hspace{3mm} \hspace{3mm} 
$\mathbf{I}[m,k] = \argmax (\mathbf{p}) + k - 1$

\textbf{return} $\mathbf{D}$, $\mathbf{I}$\\
\bottomrule
\end{tabular}
\end{table}


\subsection{Multi-grid search}
\label{sec:multi_grid_search_section}

An initial grid search over a large number of grid points is time consuming. We propose a multi-grid strategy, where we first use a coarse grid, and then refine locally around each pair of knots using a finer grid.

Suppose we have found 
$\blambda^\star = \{ {\lambda}_2^\star, \dots, {\lambda}_{K}^\star \}$ as optimal solution
for the knots on a grid of size $M-1$.
We then consider, for each $k$, a local grid of size $L$ spanning the interval $[\lambda_k^*-r, \lambda_k^*+r]$; this local grid is denoted as a vector $\btheta_k$. We then optimize the log-likelihood over all the combinations of $K-1$ values, one from each $\btheta_{k}$. (The left and right end points can be taken as the minimum and maximum values of the data points, respectively.) The parameters $L$ and $r$ can of course be set by the user, although the parameter $r$ is preferably chosen so that $2r < \delta^* = \min_k (\lambda_{k+1}^* - \lambda_{k}^*)$, so that the search intervals do not overlap.


\section{Selecting the number of modes}
\label{sec:choosing_number}

Our dynamic programming approach described in the previous section relies on knowing the value of $K$. In this section, we introduce two data-driven ways for selecting $K$ when it is unknown.

\subsection{Literature on testing multi-modality}
\label{sec:testing multi-modality}

There is a significant amount of literature on testing for multi-modality. For instance, \cite{silverman1981using} introduced a way of using kernel density estimates to investigate multi-modality, where the smoothing parameter is chosen automatically. This method was followed up by \cite{mammen1992some, hall2001calibration}, who considered its calibration and asymptotic characteristics. \cite{hartigan1985dip} proposed the dip test, which measures multi-modality in a sample by the maximum difference between the empirical distribution function and a chosen uni-modal distribution. Later, \cite{muller1991excess} considered the excess mass functional, which measures excessive empirical mass in comparison with multiples of uniform distributions. From the perspective of hierarchical clustering, \cite{hartigan1992runt} designed a test for uni-modality based on the runt size in single linkage hierarchical clustering, defined as the number of points in a cluster's smallest subcluster. 
For an early survey see, e.g., \citep{fischer1994testing}.
%

Testing for the significance of a mode was considered by \cite{chacon2013data, duong2008feature, genovese2016non, godtliebsen2002significance, rufibach2010block, minnotte1997nonparametric, cheng1999mode}, among others. \cite{hengartner1995finite} developed a lower confidence bound for the number of modes.
Multi-scale methods aiming at partitioning the real line into modal intervals, such as those of \cite{minnotte1993mode}, \cite{chaudhuri1999sizer, chaudhuri2000scale}, and \cite{dumbgen2008multiscale}, already mentioned in the introduction, are based on testing the significance of modes. Recently, \cite{ameijeiras2019mode} concluded that none of the existing methods for selecting the number of modes that are based on testing provides satisfactory performance in practice; they proposed a new method based on a combination of smoothing and excess mass.

Closely related to the problem of selecting the number of modes is the problem of selecting the number of clusters \citep{rousseeuw1987silhouettes, wang2010consistent} or the number of components in a mixture model \citep{mclachlan1987bootstrapping, mclachlan2014number, celisse2014optimal}.

\subsection{Our approaches}

Any of these existing methods for choosing the number of modes can be used to select the number of modal intervals $K$. Nonetheless, we offer two options based on our dynamic programming approach for fitting a multi-modal density. 

\subsubsection{Measure of fit}
\label{sec:measure_of_fit_subsection}

The first approach is in the spirit of the dip test \citep{hartigan1985dip}.

Let $\Delta(F,G)$ be a measure of fit for comparing the distribution functions $F$ and $G$. A famous example is the supnorm, i.e., $\Delta(F,G) = \sup_x |F(x) - G(x)|$, which is used in the Kolmogorov--Smirnov test. We find it particularly easy to intuit. Let $\hat F$ denote the empirical distribution of $\mathbf{x}$ and $\hat F_K$ denote the distribution function returned by our method. The general idea is to ``look'' at how $\Delta(\hat F_K, \hat F)$ varies with $K$. We opt for a simple approach: we set a threshold $\tau$, for example, $\tau = 0.01$, and select 
\[K^* = \min\{K : \Delta(\hat F_K, \hat F) \le \tau\}.\]

\subsubsection{Cross-validation}
\label{sec:cross_validation_subsection}
The second approach is based on cross-validation procedure. For a comprehensive review cross-validation procedures used in model selection, see the review paper by \cite{arlot2010survey}. 

We describe $R$-fold cross-validation. In that process, in each round, we fit a $K$-modal density on the data points in the folds playing the role of training, and then evaluate the log-likelihood of the fitted density on the data points in the validation fold. We then take an average of these log-likelihood values over all $R$ rounds. (We note that the log-likelihood could also be replaced by the integrated mean squared error.)

\begin{remark}
This method is, obviously, more computationally intensive.
Therefore, one might want to use a coarser dynamic programming search grid when selecting the number of modal intervals by the cross-validation, and then use a more refined grid when the choice of $K$ has been settled.
\end{remark}

\section{Numerical experiments}
\label{sec:k_mode_experiment}

In this section, we provide some numerical experiments to illustrate our dynamic programming method and its application to determining the number of modal intervals.

\subsection{Application to a real dataset}

We first apply our method to the Old Faithful Geyser dataset \citep{azzalini1990look} in the {\sf MASS} package of {\sf R}. This dataset contains $272$ waiting times, in minutes, indicating the duration between eruptions for the Old Faithful Geyser in Yellowstone National Park, Wyoming, USA. For this dataset, we consider fitting a density with the widely used choice of $K=2$. We use a search grid of size $M=5K$. We also use the multi-grid search in \secref{multi_grid_search_section}, with $L=15$ and $r = \delta^*(1/2 - 1/2L)$. 

In addition to our method, we provide the estimates of the regular kernel density estimator(KDE) with bandwidth chosen by the ``rule of thumb'' as in \citep[page 28]{silverman1986density}. We also provide the density estimate returned by the method of \citep{wolters2018enforcing}, as  implemented in the {\sf R} package {\sf scdensity} \citep{wolters2018practical}. 

\begin{figure}[!h]
\centering
\includegraphics[scale=0.55]{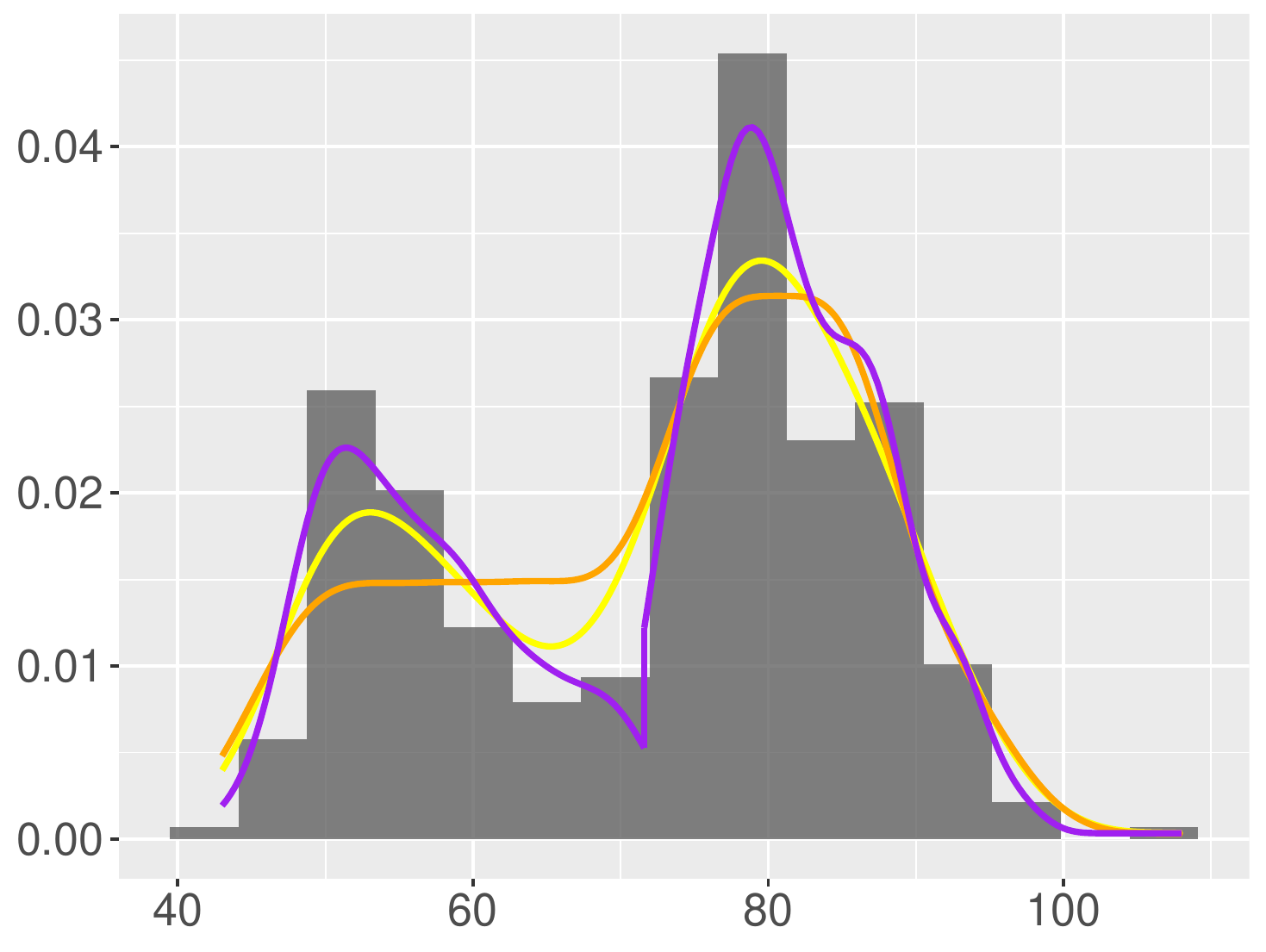}
\includegraphics[scale=0.58]{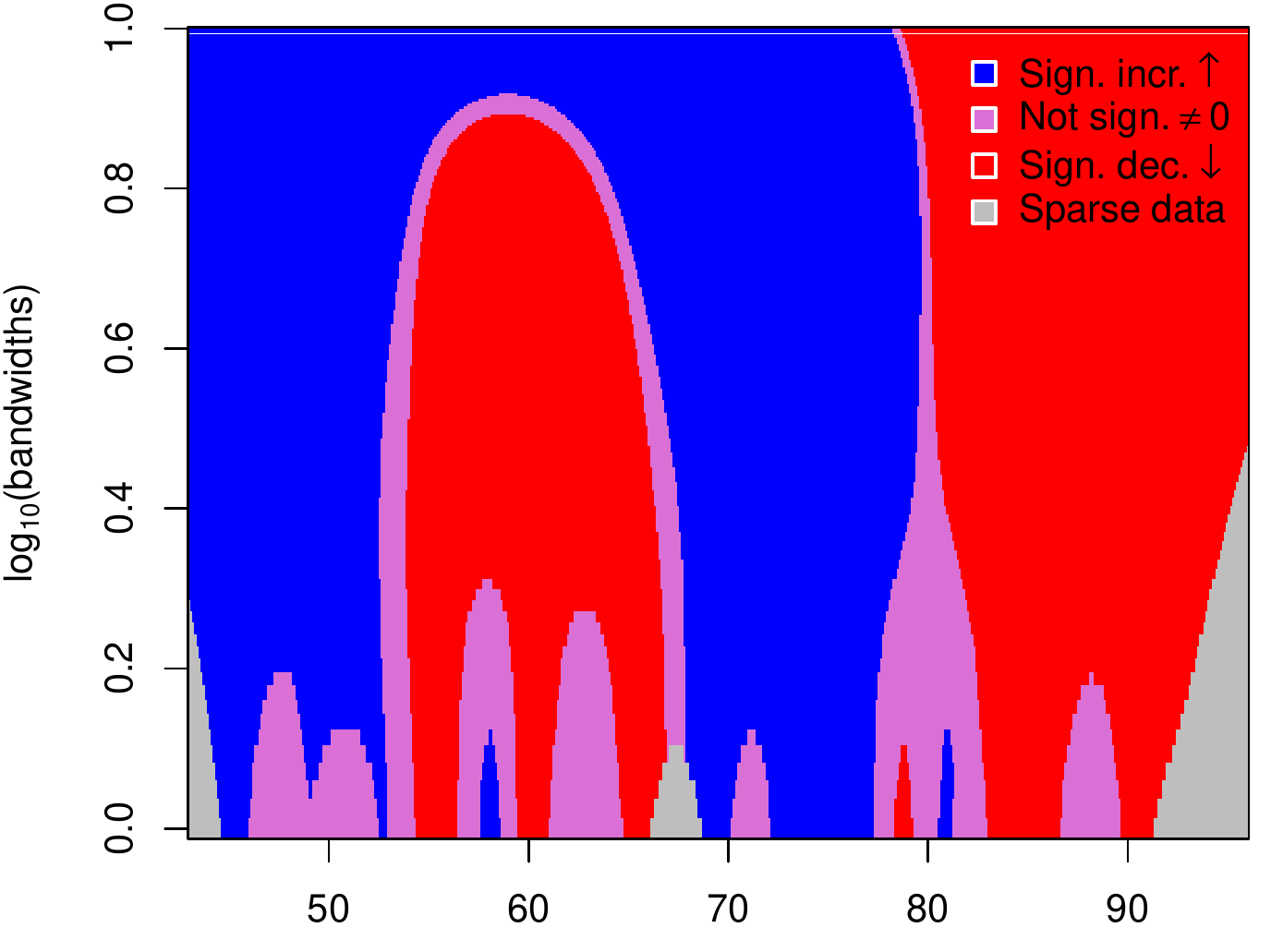}
\caption[Geyser dataset with different methods for fitting densities.]{
Left: Histogram of waiting times, in minutes, of the Geyser dataset, overlaid with density estimators acquired by regular KDE (in yellow), by additive curve (in orange), and by dynamic programming (in purple).
Right: SiZer map of the Geyser dataset.
}
\label{fig:geyser_dataset}
\end{figure}

The histogram of the waiting times, and the $3$ resulting densities, are shown in \figref{geyser_dataset}. In addition, we also provide the SiZer map \citep{chaudhuri1999sizer}, acquired from the \textsf{multimode} package \citep{ameijeiras2021multimode}. The regular KDE achieves a negative log-likelihood of $\approx 1160$; the method of \cite{wolters2018enforcing} achieves a negative log-likelihood of $\approx 1167$; and our method achieves a negative log-likelihood of $\approx 944$. We also note that in this dataset, our method gives a density that seem to fit the data better, both in high and low density regions.

\subsection{Selecting the number of modal intervals}

We next consider selecting the number of modal intervals when it is not known. To do so, we generate random mixture models, and check the performance of our method on these random mixtures.

We start by considering random Gaussian mixture models. We first uniformly sample the number of components from $\{1,2,\dots,5\}$. Once given the number of components, we randomly generate the centers uniformly at random from $[0,10]$, and we generate the standard deviations from the exponential distribution with rate $1$. We note that for each of these randomly generated Gaussian mixture models, the actual number of modes is not always the same as the number of components, and therefore we also compute the number of modes in each model by the number of locations where the derivative changes from increasing to decreasing (computed numerically on a fine grid). For each randomly generated density, we generate a sample of size $n = 10^4$. We consider here possible candidates of $K$ in $\{1,2,\dots,5\}$. (This is mainly for computational considerations, as our experiments necessarily involve replicates. In practice, the range of $K$ can be much larger.) For each possible value of $K$, we use a grid of size $M=5K$. We also use the multi-grid search as in \secref{multi_grid_search_section} with $L=5$ and $r = \delta^*(1/2 - 1/2L)$. 
We generate $B = 100$ random Gaussian mixtures in total.

We then consider random Laplace mixture models. We use the same sampling method as the Gaussian case to generate the number of components, the centers, and the standard deviations of each component. We also compute the actual number of modes like in the Gaussian case. 
We generate $B = 100$ random Laplace mixtures in total. 

\begin{remark}
While ideally we would compare with other methods, the methods we are aware of are computationally so demanding (as implemented in {\sf R}) that we were not able to include them in these simulations. This includes the recent, and promising, approach of \cite{ameijeiras2019mode}.
\end{remark}

\paragraph{Measure of fit}
The variant of \secref{measure_of_fit_subsection}, based on a measure of fit, is evaluated in \tabref{number_of_modes_ks}, where we report on the number of times the number of modal intervals was correctly identified. 
(We use the supnorm as measure of fit.)

\begin{table}[!h]
\centering\small
\caption{Number of times our method, based on the measure of fit criterion, chose the correct number of modal intervals. We present the result for $\tau = 0.01$ and $\tau = 0.05$, and for the value of $\tau$ resulting in the largest number of correct selections.}
\label{tab:number_of_modes_ks}

\bigskip
\setlength{\tabcolsep}{0.2in}
\begin{tabular}{lll}

\toprule
{\bf Mixture} & 
{\bf Threshold $\tau$} & 
{\bf Number of correct selections (out of $100$)}  \\ 
\midrule
Gaussian  & $0.01$ & $76$ \\ 
  & $0.05$ & $59$ \\
  & $0.0162$ & $78$ \\
\midrule   
Laplace  & $0.01$ & $66$ \\ 
  & $0.05$ & $50$ \\
  & $0.0104$ & $68$ \\
\bottomrule
\end{tabular}
\end{table}

\paragraph{Cross-validation}
The variant of \secref{cross_validation_subsection}, based on cross-validating the log-likelihood, is evaluated in \tabref{number_of_modes_cv}, where we report on the number of times the number of modal intervals was correctly identified. We use $R=5$ folds, and set $L=3$ for the size of the subgrids. We adopt a greedy (forward selection) strategy where, at any given value $K$, if increasing the number of modal intervals from $K$ to $K+1$ does not result in an improvement of the log-likelihood by more than 1\%, we stop and return the current value of $K$.

\begin{table}[!h]
\centering\small
\caption{Number of times our method, based on the measure of fit criterion, chose the correct number of modal intervals. We consider stopping when the improvement in the log-likelihood is less than 1\%, and also give the improvement percentage value that results in the largest number of correct selections.}
\label{tab:number_of_modes_cv}

\bigskip
\setlength{\tabcolsep}{0.2in}
\begin{tabular}{lll}

\toprule
{\bf Mixture} & 
{\bf Improvement percentage} & 
{\bf Number of correct selections (out of $100$)}  \\ 
\midrule
Gaussian  & 1\% & $66$ \\ 
  & 0.08\% & $75$ \\
\midrule   
Laplace  & 1\% & $53$ \\ 
  & 0.17\% & $62$ \\
\bottomrule
\end{tabular}
\end{table}

\section{Discussion}
\label{sec:conclusion_kmode}

In this paper, we considered a dynamic programming approach to fitting a probability density function when it is constrained to have a given number of modal intervals. Based on this approach, we then gave two data-driven ways for selecting the number of modal intervals when it is not given. 

It is in principle possible to develop a similar dynamic programming approach to fitting a density whose shape is constrained in a certain way over a number of intervals defining a partition of the real line. For example, we can consider fitting a density which is log-concave \citep{walther2009inference} in each of $K$ intervals defining a partition, where it would provide an alternative to log-concave mixture modeling \citep{chang2007clustering, pu2020algorithm}.
Importantly, though, this dynamic programming approach seems limited to real-valued data.


\small
\bibliographystyle{chicago}
\bibliography{multi_mode.bib}

\end{document}